\documentclass[12pt]{article}
\usepackage{amsmath}
\usepackage{natbib}

\begin{document}

\title{Inference on $P(Y<X)$ in Bivariate Rayleigh Distribution}
\author{$^*$A. Pak , $^\dag$N. B. Khoolenjani , $^\ddag$A. A. Jafari \\
\noindent {\small $^{*}$Department of Statistics, Shahid Chamran University, Ahvaz, Iran.}\\
{\small $^{\dag}$Department of Statistics, Shiraz University, Shiraz, Iran.}\\
{\small $^{\ddag}$Department of Statistics, Yazd University, Yazd, Iran.}\\
}
\maketitle

\begin{abstract}
This paper deals with the estimation of reliability $R=P(Y<X)$ when $X$ is a random strength of a component subjected to a random stress $Y$ and $(X,Y)$ follows a bivariate Rayleigh distribution. The maximum likelihood estimator of $R$ and its asymptotic distribution are obtained. An asymptotic confidence interval of $R$ is constructed using the asymptotic distribution. Also, two confidence intervals are proposed based on Bootstrap method and a computational approach. Testing of the reliability based on asymptotic distribution of $R$ is discussed. Simulation study to investigate performance of the confidence intervals and tests has been carried out. Also, a numerical example is given to illustrate the proposed approaches.
\end{abstract}


\bigskip

\section{Introduction}

\label{Introduction} \citet{Rayleigh} observed  that the sea waves follow no law because of the complexities of the sea, but it has been seen that the probability distributions of wave heights, wave length, wave induce pitch, wave and heave motions of the ships follow the Rayleigh distribution. Also in
reliability theory and life testing experiments, the Rayleigh distribution plays an important role.

A random variable $X$ is said to have the Rayleigh distribution with the scale parameter $\theta$ and it will be denoted by $RA(\theta)$, if its probability density function (PDF) is given by
\begin{equation}
f(x;\theta )=2\theta xe^{-\theta x^{2}},\ \ \ \ \ \theta>0, \ \ x>0.
\end{equation}

The cumulative distribution function  and survival function corresponding to (1) for $x>0$, respectively, are
\begin{equation}
F(x;\theta )=1-e^{-\theta x^{2}} ,~\ \ \ \ \ S(x;\theta)=e^{-\theta x^{2}}.
\end{equation}

Suppose  $U_{0}$ follows  $(\sim )$ $RA(\lambda _{0})$, $U_{1}\sim RA(\lambda _{1})$, $U_{2}\sim RA(\lambda _{2})$  and they are independent.  Define $X=\min \{U_{0},U_{1}\}$  and  $Y=\min\{U_{0},U_{2}\}$. Then the bivariate vector $(X,Y)$ has the bivariate Rayleigh  (BVR) distribution with the parameters $\lambda_{0}$, $\lambda _{1}$ and $\lambda _{2}$ and it will be denoted by $BVR(\lambda _{0},\lambda _{1},\lambda _{2})$.

If $(X,Y)\sim BVR(\lambda _{0},\lambda _{1},\lambda _{2})$ then their joint survival function takes the following form
\begin{eqnarray}
\bar{F}_{(X,Y)}(x,y) =P(X>x,Y>y)=e^{-\lambda _{1}x^{2}-\lambda_{2}y^{2}-\lambda _{0}(\max (x,y))^{2}}.
\end{eqnarray}
The random variables $X$ and $Y~$are independent iff $\lambda_{0}=0$. The marginals of the random variables $X$and $Y$are Rayleigh with parameters
 $\lambda _{1}+\lambda _{0}~$and $\lambda _{2}+\lambda _{0}$, respectively. The survival function of $\min (X,Y)$ is obtained by
\begin{eqnarray}
P(\min (X,Y)>x)=P(X>x,Y>x) =\bar{F}_{(X,Y)~}(x,x)=e^{-(\lambda_{0}+\lambda _{1}+\lambda _{2}) x^{2}},
\end{eqnarray}
which is the survival function of Rayleigh with parameter $\lambda _{0}+\lambda _{1}+\lambda _{2}$.

In stress-strength model, the stress $(Y)$ and the strength $(X)$ are treated as random variables and the reliability of a component during a given period is taken to be the probability that its strength exceeds the stress during the entire interval. Due to the practical point of view of reliability of stress-strength model, the estimation problem of $R=P(Y<X)$ has attracted the attention of many authors. \citet{Church}, \citet{Downtown}, \citet{Govidarajulu}, \citet{Woodward} and \citet{Owen} considered the estimation of $R$ when $X$ and $Y$ are normally distributed. \citet{Tong} considered the problem of estimating $R$, when $X$ and $Y$ are independent exponential random variables. \citet{Awad} determined the maximum likelihood estimator (MLE) of $R$ when $X$ and $Y$ have bivariate exponential
distribution. \citet{Constantine} considered the estimation of $R$ when $X$ and $Y$ are independent gamma random variables. \citet{Ahmad} and \citet{SurlesPadgett} considered the estimation of $R$ when $X$ and $Y$ are Burr type $X$ random variables. The theoretical and practical results on the theory and applications of the stress-strength relationships are collected in \citet{Kotz}. Estimation of $P(Y<X)$ from logistic (Nadarajah, 2004a), Laplace (Nadarajah, 2004b), beta (Nadarajah, 2005a), and gamma (Nadarajah, 2005b) distributions are also studied. \citet{Kundu} considered the estimation of $R$ when $X$ and $Y$ have generalized exponential distribution. Inferences on reliability in two-parameter exponential stress-strength model (Krishnamoorthy et al., 2007) and ML estimation of system reliability for Gompertz distribution (Sara\c{c}o$\check{\mbox{g}}$lu and Kaya, 2007) are considered. \citet{Kakade} considered the estimation of $R$ for exponentiated Gumbel distribution. \citet{Rezaei} considered the estimation of $R$ when $X$ and $Y$ are two independent generalized Pareto distributions.

In stress-strength analysis, usually, it is assumed that $X$ and $Y$ are independent. But any study on twins or on failure data recorded twice on  the same system naturally leads to bivariate data. For example, Meintanis (2007) considered soccer data from UEFA Champions League for the years 2004-2005 and 2005-2006 and also studied the white blood cells  counts of 60 patients. There are  not many publications on inference about stress-strength model when $X$ and $Y$ are dependent. Hanagal (1996) derived the estimator of $R=P(Y<\max(X_{1},X_{2}))$ for the case that $X_{1}$ and $X_{2}$ are bivariate exponentially distributed and stochastically independent with $Y$ that follows an exponential distribution.
 Hanagal (1997) found  estimating reliability of a component based on maximum likelihood estimators for a bivariate Pareto distribution.
 Nadarajah and Kotz (2006) studied the estimation of $P(Y<X)$ from bivariate exponential distributions.

The main aim of this paper is to discuss the inference of $R=P(Y<X)$ when $X$ is a random strength of a component subjected to a random stress $Y$ and $(X,Y)$ follows BVR distribution. In Section 2, the MLE of reliability $R$ is obtained. The asymptotic distribution of the MLE of $R$ is given and different confidence intervals are proposed in Section 3. Testing of the reliability based on a step by step computational approach is provided in Section 4. The different proposed methods are compared using Monte Carlo simulations and the results are reported in Section 5. Also, a numerical example is given to illustrate the proposed approaches.

\section{ MLE of $R$}
\label{Maximum} Suppose $Y$ and $X$ represent the random variables of stress and strength of a component, respectively, and $(X,Y)$ follows BVR distribution with
survival function given by (3). Then it can be easily seen that
\begin{eqnarray}\label{eqR}
R =P(Y<X)=\frac{\lambda _{2}}{\lambda _{0}+\lambda _{1}+\lambda_{2}}.
\end{eqnarray}
Let $(X_{i},Y_{i}),~i=1,...,n$, be a random sample of size $n$ from BVR distribution and $S_{1}$ be the random number of observations with $ y_{i}<x_{i}$ in the sample of size $n$. Then the distribution of $S_{1}$ is binomial $(n,R)$. The natural estimate of $R~$ is which is given by
\begin{equation}
\tilde{R}=\frac{1}{n}S_{1},
\end{equation}
which has the asymptotic distribution $N(R,\frac{R(1-R)}{n})$.

The MLE $\hat{R}$ of $R$ is
given as follows:
\begin{equation}
\hat{R}=\frac{\hat{\lambda}_{2}}{\hat{\lambda}_{0}+\hat{\lambda}_{1}+\hat{\lambda}_{2}}.
\end{equation}

For obtaining an explicit formula for $\hat{R}$, it is necessary
to determine the MLE's of $\lambda _{0},\lambda _{1}$ and $\lambda
_{2}$. Let $(x_{i},y_{i})$, $ i=1,...,n$, be the observations based
on a random sample of size $n$. Also, let

\begin{itemize}
\item $n_{0}=$ number of observations with $x_{i}=y_{i}$.

\item $n_{1}=$ number of observations with $x_{i}<y_{i}$.

\item $n_{2}=$ number of observations with $y_{i}<x_{i}$.
\end{itemize}
\bigskip

\noindent Then the log-likelihood function of the observed sample
is given by:
\begin{eqnarray}
L^{\ast } &=&(2n-n_{0})\log 2+\sum\limits_{j=0}^{2}n_{i}\log \lambda _{i}+n_{1}\log (\lambda _{2}+\lambda _{0})+n_{2}\log(\lambda _{1}+\lambda _{0})\nonumber \\
&&+\sum\limits_{i=1}^{n}(\log x_{i}+\log y_{i})-\sum\limits_{i\in S}\log (\max (x_{i},y_{i}))  \nonumber \\
&&-\lambda _{1}\sum\limits_{i=1}^{n}x_{i}^{2}-\lambda _{2}\sum\limits_{i=1}^{n}y_{i}^{2}-\lambda _{0}\sum\limits_{i=1}^{n}(\max (x_{i},y_{i}))^{2}.
\end{eqnarray}

The first moments or expectations of the random numbers $N_{0},N_{1}$ and $N_{2}$, are as follows:
\begin{eqnarray}
E[N_{0}]&=&(1-(\phi _{1}+\phi _{2}))n,\\
E[N_{1}]&=&\frac{\lambda _{1}}{\lambda _{1}+\lambda _{0}}(1-\phi_{1})n,\\
E[N_{2}]&=&\frac{\lambda _{2}}{\lambda _{2}+\lambda _{0}}(1-\phi_{2})n,
\end{eqnarray}
where
\begin{eqnarray}
\phi _{1}&=&\frac{\lambda _{2}}{\lambda _{1}+\lambda _{2}+\lambda_{0}},\\
\phi _{2}&=&\frac{\lambda _{1}}{\lambda _{1}+\lambda _{2}+\lambda_{0}}.
\end{eqnarray}

The MLE's of $\lambda _{0},\lambda _{1}$ and $\lambda _{2}$, say $\hat{\lambda}_{0},\hat{\lambda}_{1}$ and $\hat{\lambda}_{2}$ respectively, can be obtained as the solutions of the following system of equations:
\begin{eqnarray}
&&\frac{n_{0}}{\lambda _{0}}+\frac{n_{2}}{\lambda _{1}+\lambda _{0}}+\frac{n_{1}}{\lambda _{2}+\lambda _{0}}-\sum\limits_{i=1}^{n}(\max(x_{i},y_{i}))^{2}=0,\\
&&\frac{n_{1}}{\lambda _{1}}+\frac{n_{2}}{\lambda _{1}+\lambda _{0}}-\sum\limits_{i=1}^{n}x_{i}^{2}=0,\\
&&\frac{n_{2}}{\lambda _{2}}+\frac{n_{1}}{\lambda _{2}+\lambda _{0}}-\sum\limits_{i=1}^{n}y_{i}^{2}=0.
\end{eqnarray}

The above system of equations can be solved numerically either by using a Newton-Raphson procedure or by Fisher's method of scoring to obtain the MLE's
$(\hat{\lambda}_{0},\hat{\lambda}_{1},\hat{\lambda}_{2})$.

\section{Asymptotic distribution and Confidence Intervals}
\label{Asymptotic} In this section, first we obtain the asymptotic distribution of $\mathbf{\hat{\lambda}}=(\hat{\lambda}_{0},\hat{\lambda}_{1},\hat{\lambda}_{2})$ and then we derive the asymptotic distribution of $\hat{R}$.
Based on the asymptotic distribution of $\hat{R}$, we obtain the asymptotic confidence interval of $R$. Let us denote the Fisher information matrix of $\mathbf{\lambda } =(\lambda _{0},\lambda _{1},\lambda _{2})$ as $\mathbf{I(\mathbf{\lambda })}$. Therefore,
\begin{equation}
\mathbf{I(\lambda )}=- \left[
\begin{array}{ccc}
 E\left( \frac{\partial ^{2}L^{\ast }}{\partial \lambda_{0}^{2}}\right) \ \ \ & E \left( \frac{\partial^{2} L^{\ast }}{\partial \lambda _{0}\partial \lambda _{1}}\right) \ & E\left( \frac{\partial ^{2}L^{\ast }}{\partial \lambda _{0}\partial \lambda _{2}}\right) \\
E\left( \frac{\partial ^{2}L^{\ast }}{\partial \lambda _{1}\partial \lambda _{0}}\right)\ & E\left( \frac{\partial ^{2}L^{\ast }}{\partial \lambda _{1}^{2}}\right) \ \ \ & E\left(\frac{\partial ^{2}L^{\ast }}{\partial \lambda _{1}\partial \lambda _{2}}\right) \\
E\left( \frac{\partial ^{2}L^{\ast }}{\partial \lambda _{2}\partial \lambda _{0}}\right) \ & E\left( \frac{\partial ^{2}L^{\ast }}{\partial \lambda
_{2}\partial \lambda _{1}}\right)\ & E\left( \frac{\partial ^{2}L^{\ast }}{\partial \lambda _{2}^{2}}\right)
\end{array}
\right].
\end{equation}
Moreover,
\begin{eqnarray*}
&&E\left( \frac{\partial ^{2}L^{\ast }}{\partial \lambda _{0}^{2}}\right) =\frac{n(\phi _{1}+\phi _{2})-n}{\lambda _{0}^{2}}-\sum\limits_{i=1}^{2}\frac{
n\phi _{i}}{(\lambda _{i}+\lambda _{0})^{2}}, \nonumber \\
&&E\left( \frac{\partial ^{2}L^{\ast }}{\partial \lambda _{1}^{2}}\right) =%
\frac{n\phi _{1}-n}{\lambda _{1}(\lambda _{1}+\lambda _{0})}-\frac{n\phi _{1}}{(\lambda _{1}+\lambda _{0})^{2}}, \nonumber \\
&&E\left( \frac{\partial ^{2}L^{\ast }}{\partial \lambda _{2}^{2}}\right) =\frac{n\phi _{2}-n}{\lambda _{2}(\lambda _{2}+\lambda _{0})}-\frac{n\phi _{2}}{%
(\lambda _{2}+\lambda _{0})^{2}},  \nonumber \\
&&E\left( \frac{\partial ^{2}L^{\ast }}{\partial \lambda_{1}\partial \lambda _{0}}\right) =E\left( \frac{\partial ^{2}L^{\ast }}{\partial \lambda _{0}\partial \lambda _{1}}\right)=-\frac{n\phi _{1}}{(\lambda _{1}+\lambda_{0})^{2}},  \nonumber \\
&&E\left( \frac{\partial ^{2}L^{\ast }}{\partial \lambda_{2}\partial \lambda _{0}}\right) =E\left( \frac{\partial^{2}L^{\ast }}{\partial \lambda _{0}\partial \lambda _{2}}\right) =-\frac{n\phi _{2}}{(\lambda _{2}+\lambda _{0})^{2}}, \nonumber \\
&&E\left( \frac{\partial ^{2}L^{\ast }}{\partial \lambda _{1}\partial \lambda _{2}}\right) =E\left( \frac{\partial ^{2}L^{\ast }}{\partial \lambda _{2}\partial \lambda _{1}}\right)=0. \nonumber
\end{eqnarray*}

The above Fisher information matrix is positive definite and by the asymptotic results for the MLE, we arrive at the following theorem:

\bigskip

\noindent{\bf Theorem 1.} {\it As $n\rightarrow \infty$, then
\begin{equation}
\sqrt{n}(\mathbf{\hat{\lambda}-\lambda }){\rightarrow }N_{3}(\mathbf{0},\mathbf{J^{-1}(\lambda )}),
\end{equation}
where $\mathbf{J(\lambda)}=n^{-1}\mathbf{I(\lambda )}$ and  $\mathbf{J^{-1}(\lambda )}$ is the inverse of $\mathbf{J(\lambda)}$.}

\noindent {\bf Proof.} {\it The result follows straightforward from the asymptotic properties of MLE's under regularity conditions and the multivariate central limit theorem.}

\bigskip

\noindent{\bf Theorem 2.} {\it The asymptotic distribution of $\hat{R}$ is normal with the value of the first moment $R$ and the value of the variance $\Sigma$ that is given by
\begin{equation}
\Sigma =\mathbf{B^{T}}\mathbf{G}\mathbf{B},
\end{equation}
with
\begin{equation}
\mathbf{B^{T}}=\left( \frac{\partial R}{\partial \lambda _{0}},\frac{\partial R}{\partial \lambda _{1}},\frac{\partial R}{\partial \lambda_{2}}\right),
\end{equation}
and $\mathbf{G}$ is the inverse of the variance-covariance matrix of $\left(\hat{\lambda}_{0},\hat{\lambda}_{1},\hat{\lambda}_{2}\right)$.}

\noindent{\bf Proof.} {\it The result follows from invariance property of consistent asymptotically normal estimators under continuous transformation.(See \citet{Ferguson}, Section 7).}

\bigskip
\noindent{\bf Remark 1.} {\it By means of Theorem 2, an asymptotic confidence interval of $R$~ is obtained as follows:
\begin{equation}
\left( \hat{R}-Z_{1-\frac{\alpha }{2}}\sqrt{\hat{\Sigma}}~,~\hat{R}+Z_{1-\frac{\alpha }{2}}\sqrt{\hat{\Sigma}}\right).
\end{equation}}

\bigskip

\noindent{\bf Remark 2.} {\it The value $\Sigma $ of the variance can be estimated by means of the empirical Fisher information matrix and the MLE's of $\lambda _{0},\lambda _{1}$ and $\lambda _{2}$.}

\bigskip

\noindent{\bf Remark 3.} {\it Using Theorem 2, an asymptotic test of size $\alpha $ rejects the null hypothesis $H_{0}:R=R_{0}$ against $H_{1}: R>R_{0}$, if
\begin{equation}
(\hat{R}-R_{0})>Z_{1-\alpha }\sqrt{\hat{\Sigma}}
\end{equation}
where $Z_{1-\alpha }$ is the $(1-\alpha )^{th}$ quantile of the standard normal distribution. We can also obtain asymptotic tests of the desired size for alternatives $H_{1}:R<R_{0}$ or $H_{1}:R\neq R_{0}$.}

\bigskip

It is observed that the asymptotic confidence intervals do not perform very well. Therefore, we propose the following
bootstrap confidence interval.

\subsection{Bootstrap confidence interval}

In this subsection, we propose a percentile bootstrap method (Efron, 1982) for constructing confidence interval of $R$ which is as follows.

\noindent \textbf{step 1.} Generate random sample $\left(
x_{1},y_{1}\right) ,...,\left( x_{n},y_{n}\right)$ from
$BVR\left( \lambda _{0},\lambda _{1},\lambda _{2}\right) $ and
compute $\hat{\lambda}_{0}$, $\hat{\lambda}_{1}$ and
$\hat{\lambda}_{2}$.

\noindent \textbf{step 2.} Using $\hat{\lambda}_{0}$, $\hat{\lambda}_{1}$ and $\hat{ \lambda}_{2}$ generate a bootstrap sample $\left( x_{1}^{\ast },y_{1}^{\ast }\right) ,...,\left(x_{n}^{\ast },y_{n}^{\ast }\right) $ from $BVR(\hat{\lambda}_{0},\hat{\lambda}_{1},$ $\hat{\lambda}_{2})$.
Based on this bootstrap sample compute bootstrap estimate of $R$
using (7), say $\hat{R}^{\ast } $.

\noindent \textbf{step 3.} Repeat step 2, NBOOT times.

\noindent \textbf{step 4.} Let $H(x)=P(\hat{R}^{\ast }\leq x)$, be
the cumulative
distribution function of $\hat{R}^{\ast }$. Define $\hat{R}_{Boot-p}(x)=H^{-1}(x)~$%
for a given $x$. The approximate $100(1-\alpha )\%~$bootstrap
confidence interval of $R$ is given by
\begin{equation}
\left( \hat{R}_{Boot-p}(\frac{\alpha
}{2})~,\hat{R}_{Boot-p}(1-\frac{\alpha }{2})~\right) .
\end{equation}

\section{ Hypothesis testing and Interval Estimation Based on a Computational Approach}
\label{Hypothesis} In this section, we use the idea of \citet{Pal} to testing the reliability and constructing confidence interval of $R$ based on the MLE. The proposed computational approach test (CAT) based on simulation and numerical computations uses the ML estimate(s), but does not require any asymptotic distribution.

\subsection{ Hypothesis Testing and the Computational Approach Test}

Suppose $(X_{1},Y_{1}),...,(X_{n},Y_{n})$ are $iid$ random samples from $ BVR(\lambda _{0},\lambda _{1},\lambda _{2})$. Our goal is to test $ H_{0}:R=R_{0}$\ \ against a suitable $H_{1}(R<R_{0}$ or $R>R_{0}$ or $ R\neq R_{0})$ at level $\alpha $. Under $H_{0}$, the log-likelihood function of the sample of size
$n$ can be expressed as
\begin{eqnarray}
L^{\ast } &=&(2n-n_{0})\log 2+n_{0}\log \lambda _{0}+n_{1}\log \lambda _{1}+n_{2}\log (\frac{R_{0}}{1-R_{0}}(\lambda _{1}+\lambda _{0}))  \nonumber \\
&&+n_{1}\log (\frac{R_{0}}{1-R_{0}}(\lambda _{1}+\lambda _{0})+\lambda _{0})+n_{2}\log (\lambda _{1}+\lambda _{0}) \nonumber\\
&& +\sum\limits_{i=1}^{n}(\log x_{i}+\log y_{i}) -\sum\limits_{i\in S}\log (\max (x_{i},y_{i}))\nonumber \\
&&-\lambda _{1}\sum\limits_{i=1}^{n}x_{i}^{2}-\frac{R_{0}}{1-R_{0}}(\lambda _{1}+\lambda _{0})\sum\limits_{i=1}^{n}y_{i}^{2}-\lambda_{0}\sum\limits_{i=1}^{n}(\max (x_{i},y_{i}))^{2}.
\end{eqnarray}

 The MLEs of $\lambda _{0}$ and $\lambda _{1} $, under $H_{0}$,$~$can be obtained as the solutions of the following system of equations:
\begin{eqnarray}
&&\frac{n_{0}}{\lambda _{0}}+\frac{2n_{2}}{\lambda _{1}+\lambda _{0}}+\frac{n_{1}R_{0}}{\lambda _{0}+\lambda _{1}R_{0}}-\sum\limits_{i=1}^{n}(\max(x_{i},y_{i}))^{2}-\frac{R_{0}}{1-R_{0}}\sum\limits_{i=1}^{n}y_{i}^{2}=0 \, \, \\
&&\frac{n_{1}}{\lambda _{1}}+\frac{2n_{2}}{\lambda _{1}+\lambda _{0}}+\frac{n_{1}}{\lambda _{0}+\lambda _{1}R_{0}}-\sum\limits_{i=1}^{n}x_{i}^{2}-\frac{R_{0}}{1-R_{0}}\sum\limits_{i=1}^{n}y_{i}^{2}=0
\end{eqnarray}

\noindent The CAT is given through the following steps:

\bigskip

\noindent\textbf{step 1.} Obtain the MLEs  $\hat{\lambda}_{0}$, $\hat{\lambda}_{1}$ and $\hat{\lambda}_{2}$ from equations (14), (15) and (16) and compute the MLE of $R$, say $\hat{R}_{ML}$, from (7).

\bigskip
\noindent\textbf{step 2.}
\begin{itemize}
 \item[(i)]  Set $R=R_{0},$ then find the MLEs $\hat{\lambda}_{0}$ and $\hat{\lambda} _{1}$ from the original data by using (25) and
(26),  and call this as the ``restricted MLE of  $(\lambda _{0},\lambda _{1})$'', denoted by $(\hat{\lambda}_{0RML},\hat{\lambda}_{1RML})$.

\item [(ii)] Generate artificial sample$(X_{1},Y_{1}),...,(X_{n},Y_{n})$,  iid from  BVR$(\hat{\lambda}_{0RML},\hat{\lambda}_{1RML}, \frac{R_{0}}{1-R_{0}}\times$ $ (\hat{\lambda}_{1RML}+\hat{\lambda}_{0RML}))$   a large number of times (say, M times). For each of replicated samples, recalculate the MLE of $R$. Let these recalculated MLE values of $R$ be $\hat{R}_{01},\hat{R}_{02},...,\hat{R}_{0M}$.

\item [(iii)] Let $\hat{R}_{0(1)}<\hat{R}_{0(2)}<...<\hat{R}_{0(M)}$ be the ordered values of $\hat{R}_{0l}$, \ $1\leq l\leq M$.
  \end{itemize}

\noindent\textbf{step 3.}

\begin{itemize}
 \item[(i)] For testing $H_{0}$ against $H_{1}:R<R_{0}$, define $\hat{R}_{l}=\hat{R}_{0(\alpha M)}$. Reject $H_{0}$ if $\hat{R}_{ML}<\hat{R}_{l}$. Alternatively, calculate the $p$-value as: $p=$ (number of $\hat{R}_{0(l)} $'s $<\hat{R}_{ML})/M$.
\item[(ii)] For testing $H_{0}$ against $H_{1}:R>R_{0}$, define $\hat{R}_{u}=\hat{R}_{0((1-\alpha )M)}$. Reject $H_{0}$ if $\hat{R}_{ML}>\hat{R}_{u}$. Alternatively, calculate the $p$-value as: $p=$ (number of $\hat{R} _{0(l)}$'s $>\hat{R}_{ML})/M$.
 \item[(iii)] For testing $H_{0}$ against $H_{1}:R\neq R_{0}$, define $\hat{R}_{u}=\hat{R}_{0((1-\frac{\alpha }{2})M)}$ and $\hat{R}_{l}=\hat{R}_{0((\frac{\alpha }{2})M)}$. Reject $H_{0}$ if $\hat{R}_{ML}$ is either greater than $\hat{R}_{u}$ or less than $\hat{R}_{l}$.
Alternatively, calculate the $p$-value as: $p=2\min (p_{1},p_{2})$ where $p_{1}=$(number of $\hat{R}_{0(l)}$'s $<\hat{R}_{ML})/M$ and $p_{2}=$(number of $\hat{R}_{0(l)}$'s $>\hat{R}_{ML})/M$.
\end{itemize}

\subsection{ Interval Estimation}

Since we have already discussed about hypothesis testing based on our suggested CAT, we can take advantage of it for constructing confidence interval of the reliability $R$, by the following steps:

\bigskip
\noindent\textbf{step 1.}

\begin{itemize}
 \item[(i)] Take a few values $R_{0}^{1},...,R_{0}^{k}$ of $R$, preferably equally spaced. [It is suggested that these values of $R$, numbered between 8 and 12 (more than better)].
 \item[(ii)] Perform a hypothesis testing of $H_{0}:R=R_{0}^{j}$ against $H_{1}^{j}:R\neq R_{0}^{j}~~(1\leq j\leq k)$ at level $\alpha$. For each $j$, obtain the cut-off points $(\hat{R}_{L}^{j},\hat{R}_{U}^{j})$ based on the CAT described earlier. Note that the bounds $(\hat{R}_{L}^{j},\hat{R}_{U}^{j})$ are based on $\hat{R}_{ML}~$when $R=R_{0}^{j}$ is considered to be the true value of $R$, $1\leq j\leq k$.
 \end{itemize}

\noindent\textbf{step 2.} Plot the lower bounds $\hat{R}_{L}^{j}$, $1\leq j\leq k$, against $R_{0}^{j}$, and then approximate the plotted curve by a suitable smooth function, say $g_{L}(R_{0})$. Similarly plot the upper bounds $\hat{R}_{U}^{j}$, $1\leq j\leq k,$ against $R_{0}^{j}$, and then approximate the plotted curve by a suitable smooth function, say $g_{U}(R_{0})$.

\bigskip

\noindent\textbf{step 3.} Finally, solve for $R_{0}$ from the equations $g_{L}(R_{0})= \hat{R}_{ML}$ and $g_{U}(R_{0})=\hat{R}_{ML}$. The two solutions of $R_{0}$ thus obtained set the boundaries of the interval estimate of $R$ with intended confidence bound $(1-\alpha)$.

\bigskip In Section 5, we apply the above mentioned computational approach for BVR distribution.

\section{Numerical studies}
In this section we first present some simulation experiments to observe the behavior of the different methods for various sample sizes and for various values of parameters. Then, a numerical example is provided for illustrating the proposed approaches to find 95\% confidence interval for $R$. The  data set has been obtained from Meintanis (2007).

\subsection{ Simulation Study}
We evaluate  the performances of the MLEs with respect to the squared error loss function in terms of biases and mean squared errors (MSEs). We consider the sample sizes $n=5,10,15,20,25$, $50$ and the parameter values $\lambda _{1}=1$, $\lambda _{2}=1$ and $\lambda_{0}=0.5,1,1.5,2,2.5$. All the results are based on 1000 replications. From the sample, we obtain the MLE of $R$ using (7). The average biases and MSEs of the MLEs are presented in Table \ref{tab:1}, based on 1000 replications. Some of the points are quite clear from this experiment. Even for small sample sizes, the performance of the MLEs are quite satisfactory in terms of biases and MSEs. When sample size increases, the MSEs decrease. It verifies the consistency property of the MLE of $R$.

We also compute the 95\% confidence intervals and estimate average lengths and coverage percentages of asymptotic confidence intervals, bootstrap confidence intervals and the confidence intervals obtained by using the computational approach given in section 4. The results are reported in Table \ref{tab:2}.

It is observed that when the sample sizes increase, the coverage percentages of  asymptotic confidence interval and computational approach increase but they are always smaller than the confidence coefficient even for samples as large as 50. The performance of the bootstrap confidence intervals are quite well and the coverage percentages of this method are close to the confidence coefficient. In fact, it is clear that the bootstrap approach works far better than the other methods.

Through simulation study, comparison of power is made for asymptotic test and the computational approach test (CAT) given in section
4. The power is determined by generating 1000 random samples of size $n=10,15,20,25,50$. The results for the test $H_{0}:R=R_{0}$\ \ against \ $R>R_{0}$ at the significance level $\alpha =0.05$ are presented in Table 3. $P_{1}$ and $P_{2} $ are referred to as powers based
on the CAT and the asymptotic test, respectively. The following points are observed from Table \ref{tab:power}:

\begin{itemize}
\item Both tests perform well with respect to the power.
\item It is clear that the CAT is almost as good as the asymptotic test. The whole idea behind the CAT has been the assertion that not
knowing or applying the sampling distribution of the MLE of $R$ does not cause much detriment as seen in Table \ref{tab:power}.
\item Both the tests are consistent in the sense that as sample sizes increase, the power of the tests show improvement.
\end{itemize}

\subsection{Real example}
Table 4 represent the football (soccer) data for the group stage of the UEFA Champion's League for the years 2004-2005 and 2005-2006. Considered matches where there was at least one goal scored by the home team and at least one goal scored directly from a penalty kick, foul kick or any other direct kick (all of them together will be called as kick goal) by any team. Here, the time in minutes of the first kick goal scored by any team is represented by $X_1$, and the first goal of any type scored by the home team is represented by $X_2$. In this case $(X_1,X_2)$ is a bivariate continuous random vector for which all possibilities are open, for example $X_1\ <\ X_2$, or $X_1\ >\ X_2$ or $X_1=X_2$.

Meintanis (2007) used the Marshal-Olkin distribution (Marshal and Olkin, 1967) to analyze these data and Kundu and Gupta (2009) re-analyzed using bivariate generalized exponential distribution.

We fitted the Rayleigh distribution to $X_1$ and $X_2$ separately. The Kolmogorov-Smirnov distances between the fitted distribution and the empirical distribution function and the corresponding $p$-values (in brackets) for $X_1$ and $X_2$ are 0.0885 (0.9341) and 0.1897 (0.1073) respectively. Based on the $p$-values, Rayleigh distribution cannot be rejected for the marginals. We also fit the BVR distribution to the data and obtained the MLE of $R$ as $\hat{R}=0.4228$. We computed the 95\% confidence intervals of $R$ based on asymptotic confidence interval, bootstrap confidence interval, and computational approach as (0.280, 0.565), (0.276, 0.571), and (0.201, 0.637), respectively.

\section{ Conclusions}

In this paper, we have addressed the inferences of $R=P(Y<X)$ when $X$ is a random strength of a component subjected to a random stress $Y$ and $(X,Y)$ follows bivariate Rayleigh distribution. It is observed that the MLE works quite well. Based on the simulation results, we recommended to use the Bootstrap confidence interval, even, when the sample size is very small. Using the asymptotic distribution of the MLE to construct confidence intervals does not work well.  We adopt a step by step computational approach to handle statistical inferences. This approach may come handy in those cases where the sampling distributions are not easy to derive or extremely complicated. The simulation studies on powers of  the computational approach test  and the asymptotic test show that both tests perform satisfactory.

\section*{Acknowledgements}
The authors would like to thank two referees for many constructive suggestions.

\bigskip

\begin{table}[h]
\caption{Biases and MSEs of the MLEs of $R$, when $\lambda _{1}=1$, $\lambda _{2}=1$ for different values of $\lambda _{0}$.}\label{tab:1}
{\small
\begin{tabular}{cccccc}
\hline $n$  & $\lambda _{0}=0.5$ & $\lambda _{0}=1$ & $\lambda_{0}=1.5$ & $\lambda _{0}=2$ & $\lambda _{0}=2.5$\\ \hline
 5 & -0.0038(0.0378) & -0.0074(0.0351) & -0.0110(0.0249) & -0.0183(0.0202) & -0.0206(0.0189) \\
10 & -0.0067(0.0183) & -0.0085(0.0165) & -0.0143(0.0138) & -0.0173(0.0112) & -0.0185(0.0094) \\
15 & -0.0023(0.0109) & -0.0070(0.0091) & -0.0093(0.0089) & -0.0117(0.0084) & -0.0129(0.0071) \\
20 & -0.0019(0.0088) & -0.0046(0.0073) & -0.0058(0.0064) & -0.0059(0.0052) & -0.0066(0.0047) \\
25 & -0.0018(0.0075) & 0.0016(0.0062)  & -0.0044(0.0045) & -0.0048(0.0039) & -0.0051(0.0034) \\
50 & -0.0008(0.0034) & -0.0017(0.0031) & -0.0023(0.0029) & -0.0026(0.0024) & -0.0028(0.0019) \\
\hline
\end{tabular}
}
\end{table}

\begin{table}
\caption{Average confidence lengths and coverage percentages, when $\lambda _{1}=1$, $\lambda _{2}=1$ for different values of
$\lambda _{0}$.} \label{tab:2}
{\small
\begin{tabular}{ccccccc}
\hline $n$& & $\lambda _{0}=0.5$ & $\lambda _{0}=1$ &
$\lambda _{0}=1.5$ & $\lambda _{0}=2$ & $\lambda _{0}=2.5$\\
\hline
 5 &a& 0.5477(0.835) & 0.5326(0.821) & 0.5274(0.817) & 0.4875(0.803) & 0.4519(0.798) \\
   &b& 0.5360(0.942) & 0.5029(0.930) & 0.4823(0.947) & 0.4535(0.933) & 0.4487(0.922) \\
   &c& 0.5532(0.832) & 0.5628(0.811) &0.5319(0.806)  &0.4693(0.789)  & 0.4418(0.790) \\ \hline
10 &a& 0.4853(0.902) & 0.4803(0.891) & 0.4694(0.882) & 0.4532(0.879) & 0.4489(0.874) \\
   &b& 0.5097(0.929) & 0.5026(0.962) & 0.4931(0.958) & 0.4737(0.928) & 0.4515(0.934) \\
   &c& 0.4954(0.893) &0.4878(0.887)  &0.4912(0.869)  &0.4564(0.853)  & 0.4783(0.866) \\ \hline
15 &a& 0.4096(0.918) & 0.4033(0.901) & 0.3987(0.890) & 0.3829(0.886) & 0.3712(0.883) \\
   &b& 0.4086(0.948) & 0.4023(0.954) & 0.4009(0.969) & 0.3973(0.961) & 0.3971(0.951) \\
   &c& 0.4170(0.905) & 0.4219(0.890) & 0.4016(0.871) & 0.3859(0.884) & 0.3663(0.876) \\ \hline
20 &a& 0.3563(0.928) & 0.3441(0.917) & 0.3378(0.911) & 0.3215(0.908) & 0.3195(0.907) \\
   &b& 0.3527(0.957) & 0.3504(0.962) & 0.3489(0.949) & 0.3360(0.953) & 0.3312(0.937) \\
   &c& 0.3859(0.920) & 0.3519(0.903) & 0.3451(0.910) & 0.3586(0.891) & 0.3219(0.885) \\ \hline
25 &a& 0.3206(0.932) & 0.3179(0.931) & 0.3032(0.934) & 0.2918(0.928) & 0.2845(0.921) \\
   &b& 0.3274(0.964) & 0.3257(0.942) & 0.3110(0.951) & 0.2936(0.969) & 0.2905(0.935) \\
   &c& 0.3216(0.925) & 0.3163(0.909) & 0.3144(0.918) & 0.2991(0.922) & 0.2932(0.893) \\ \hline
50 &a& 0.2299(0.933) & 0.2257(0.931) & 0.2124(0.928) & 0.1983(0.925) & 0.1911(0.922) \\
   &b& 0.2296(0.937) & 0.2273(0.947) & 0.2203(0.953) & 0.2145(0.948) & 0.2118(0.950) \\
   &c& 0.2314(0.932) & 0.2468(0.917) & 0.2183(0.926) & 0.2208(0.924) & 0.1963(0.915) \\
\hline
\end{tabular}\\
(a: The first row represents the average confidence lengths based
on the asymptotic distributions of the MLEs. The corresponding
coverage percentages are reported within brackets. b and c:
Similarly the second and third rows represent the results for
Bootstrap method and the computational approach in 4.2,
respectively.)
}
\end{table}

\begin{table}
\caption{Power of the CAT and the asymptotic test, $\alpha =0.05$.} \label{tab:power}
{\small
\begin{tabular}{ccccccccccccccc}

\hline
& \multicolumn{14}{c}{$n$}\\
\cline{2-15}
  &\multicolumn{2}{c}{ 10}&& \multicolumn{2}{c}{ 15}&& \multicolumn{2}{c}{20}&& \multicolumn{2}{c}{ 25}&&\multicolumn{2}{c}{ 50} \\
 \cline{2-3} \cline{5-6} \cline{8-9} \cline{11-12} \cline{14-15}
$R$  & $P_{1}$  &$P_{2}$ & &$P_{1}$&$P_{2}$  &&$P_{1}$  &$P_{2}$&&$P_{1}$&$P_{2}$ &&$P_{1}$&$P_{2}$ \\
\hline

0.500 &0.065 & 0.070& & 0.059 & 0.069& &0.066 &  0.061 && 0.053  &0.054 && 0.048 & 0.049  \\
0.534 &0.119 & 0.122 && 0.138 & 0.142& & 0.136&  0.144& & 0.147 & 0.159 & &0.176 & 0.181  \\
0.562 &0.153 & 0.161 & &0.269  &0.176 & &0.187 & 0.195& & 0.221 & 0.229 && 0.321  &0.330  \\
0.600 &0.277 & 0.282 & &0.308 & 0.314  &&0.359  &0.362 && 0.408 & 0.413  &&0.558 & 0.561  \\
0.636 &0.330 & 0.329 & &0.410  &0.412& & 0.457&  0.468  &&0.571  &0.573  &&0.681  &0.692  \\
0.666 &0.490&  0.497 && 0.562 & 0.570 && 0.653 & 0.655& & 0.728 & 0.733 & &0.935  &0.936  \\
0.714 &0.639 & 0.644 && 0.758  &0.766  &&0.866 & 0.868 && 0.916 & 0.919 && 0.993   &0.995  \\
0.777 &0.856 & 0.858 && 0.954 & 0.956& & 0.973 & 0.977 & &0.981 & 0.983 && 0.996      &1      \\
0.833 &0.970 & 0.972 & &0.991 & 0.998 && 0.995 & 0.998 && 1     & 1     & &1     & 1      \\
0.882 &0.996 & 0.997 && 1      &1     && 1     & 1     & &1     & 1    & & 1     & 1      \\
\hline
\end{tabular}
}
\end{table}

\begin{table}
\begin{center}
\caption{UEFA Champion's League data.} \label{tab:ex}
{\small

\begin{tabular}{|l|c|c|c|l|c|c|} \hline
2005-2006 & $X_1$ & $X_2$ &  & 2004-2005 & $X_1$ & $X_2$ \\ \hline
Lyon-Real Madrid & 26 & 20 &  & Internazionale-Bremen & 34 & 34 \\
Milan-Fenerbahce & 63 & 18 &  & Real Madrid-Roma & 53 & 39 \\
Chelsea-Anderlecht & 19 & 19 &  & Man. United-Fenerbahce & 54 & 7 \\
Club Brugge-Juventus & 66 & 85 &  & Bayern-Ajax & 51 & 28 \\
Fenerbahce-PSV & 40 & 40 &  & Moscow-PSG & 76 & 64 \\
Internazionale-Range & 49 & 49 &  & Barcelona-Shakhtar & 64 & 15 \\
Panathinaikos-Bremen & 8 & 8 &  & Leverkusen-Roma & 26 & 48 \\
Ajax-Arsenal & 69 & 71 &  & Arsenal-Panathinaikos & 16 & 16 \\
Man. United-Benfica & 39 & 39 &  & Dynamo Kyiv-Real Madrid & 44 & 13 \\
Real Madrid-Rosenborg & 82 & 48 &  & Man. United-Sparta & 25 & 14 \\
Villarreal-Benfica & 72 & 72 &  & Bayern-M. Tel-Aviv & 55 & 11 \\
Juventus-Bayern & 66 & 62 &  & Bremen-Internazionale & 49 & 49 \\
Club Brugge-Rapid & 25 & 9 &  & Anderlecht-Valencia & 24 & 24 \\
Olympiacos-Lyon & 41 & 3 &  & Panathinaikos-PSV & 44 & 30 \\
Internazionale-Porto & 16 & 75 &  & Arsenal-Rosenborg & 42 & 3 \\
Schalke-PSV & 18 & 18 &  & Liverpool-Olympiacos & 27 & 47 \\
Barcelona-Bremen & 22 & 14 &  & M. Tel-Aviv-Juventus & 28 & 28 \\
Milan-Schalke & 42 & 42 &  & Bremen-Panathinaikos & 2 & 2 \\
Rapid-Juventus & 36 & 52 &  &  &  &  \\ \hline
\end{tabular}
}
\end{center}
\end{table}

\newpage

\hspace{5cm}


\newpage



\begin{thebibliography}{}

\bibitem[Ahmad et al. (1997)]{Ahmad}
Ahmad, K.E., Fakhry, M.E., Jaheen, Z.F., 1997. Empirical Bayes estimation of $P(Y<X)$ and characterizations of the Burr-type X model. {\it Journal of Statistical Planning and Inference} 64: 297-308.


\bibitem[Awad et al. (1981)]{Awad} Awad, A.M., Azzam, M.M. , Hamadan, M.A., 1981. Some inference results in $P(Y<X)$ in the bivariate exponential model.
{\it Communications in Statistics-Theory and Methods} 10: 2515-2524.

\bibitem[Church and Harris (1970)]{Church} Church, J.D., Harris, B., 1970. The estimation of reliability from stress strength relationships. {\it Technometrics} 12: 49-54.

\bibitem[Constantine et al. (1986)]{Constantine} Constantine, K., Tse, S., Karson. M., 1986. The estimation of $P(Y<X)$ in gamma case. {\it Communication in Statistics-Computations and Simulations} 15: 365-388.

\bibitem[Downtown (1973)]{Downtown} Downtown, F., 1973. The estimation of $P(Y<X)$ in the normal case. {\it Technometrics} 15: 551-558.

\bibitem[Efron (1982)]{Efron} Efron, B., 1982. {\it The Jackknife , the bootstrap and other resampling plans}. In: CBMS-NSF Regional Conference Series in
Applied Mathematics 38, SIAM, Philadelphia, PA.

\bibitem[Ferguson (1996)]{Ferguson} Ferguson, T.S., 1996. {\it A Course in Large Sample Theory}. New York: Chapman and Hall.

\bibitem[Govidarajulu (1967)]{Govidarajulu} Govidarajulu, Z., 1967. Two sided confidence limits for $P(Y<X)$ based on normal sample of $X$ and $Y$. {\it Sankhy$\tilde{a}$} B 29: 35-40.

\bibitem[Hanagal (1996)]{Hanagal} Hanagal, D.D., 1996. Estimation of system reliability from Stress-Strength relationship. {\it Communications in Statistics-Theory and Methods} 25(8): 1783-1797.

\bibitem[Hanagal (1997)]{Hanagal2} Hanagal, D.D., 1997. Note on estimation of reliability under bivariate Pareto stress-strength model. {\it Statistical Papers} 38: 453-459.

\bibitem[Kakade et al. (2008)]{Kakade} Kakade, C.S., Shirke, D.T., Kundu, D., 2008. Inference for $P(Y<X)$ in exponentiated Gumbel distribution. {\it Journal of Statistics and Applications} 3 (1-2): 121-133.

\bibitem[Kotz et al. (2003)]{Kotz} Kotz, S., Lumelskii, Y., Pensky, M., 2003. {\it The Stress-Strength model and its Generalizations: Theory and Applications}. Word Scientific, Singapore.

\bibitem[Krishnamoorthy et al. (2007)]{Krishnamoorthy} Krishnamoorthy, K., Mukherjee, S., Guo, H., 2007. Inference on reliability in two-parameter exponential stress-strength model. {\it Metrika} 65 (3): 261-273.

\bibitem[Kundu and Gupta (2005)]{Kundu} Kundu, D., Gupta, R.D., 2005. Estimation of $P(Y<X)$ for Generalized Exponential Distributions. {\it Metrika} 61(3): 291-308.

\bibitem[Kundu and Gupta (2009)]{Kundu2} Kundu, D., Gupta, R.D., 2009. Bivariate generalized exponential distribution. {\it Journal of Multivariate Analysis} 100, 581-593.

\bibitem[Marshall and Olkin (1967)]{Marshall} Marshall, A. W. and Olkin, I., 1967. A multivariate exponential distribution. {\it Journal of the American Statistical Association} 62, 30-44.

\bibitem[Meintanis (2007)]{Meintanis} Meintanis, S. G., 2007. Test of fit for Marshall-Olkin distributions with applications. {\it Journal of Statistical Planning and Inference} 137, 3954-3963.

\bibitem[Nadarajah, 2004a]{Nadarajah} Nadarajah, S., 2004a. Reliability for logistic distributions. {\it Elektronnoe Modelirovanie} 26 (3): 65-82.

\bibitem[Nadarajah, 2004b]{Nadarajah2} Nadarajah, S., 2004b. Reliability for Laplace distributions. {\it Mathematical Problems in Engineering} 2: 169-183.

\bibitem[Nadarajah, 2005a]{Nadarajah3} Nadarajah, S., 2005a. Reliability for some bivariate beta distributions. {\it Mathematical Problems in Engineering} 1: 101-111.

\bibitem[Nadarajah, 2005b]{Nadarajah4} Nadarajah, S., 2005b. Reliability for some bivariate gamma distributions. {\it Mathematical Problems in Engineering} 2: 151-163.

\bibitem[Nadarajah and Kotz, 2006]{Nadarajah5} Nadarajah, S., Kotz, S., 2006. Reliability for some bivariate exponential distributions. {\it Mathematical Problems in Engineering} 2006: 1-14.

\bibitem[Owen et al. (1977)]{Owen} Owen, D.B., Craswell, K.J., Hanson, D.L., 1977. Non-parametric upper confidence bounds for $P(Y<X)$ and confidence limits for
$P(Y<X)$ when $X$ and $Y$ are normal. {\it Journal of the American Statistical Association} 59: 906-924.

\bibitem[Pal et al. (2007)]{Pal} Pal, N., Lim, W.K., Ling, C.H., 2007. A Computational Approach to Statistical Inferences. {\it Journal of Applied Probability and Statistics} 2(1): 13-35.

\bibitem[Rayleigh (1880)]{Rayleigh} Rayleigh, L., 1880. On the resultant of a large number vibrations of the same pitch and of arbitrary phase. {\it Philosophical Magazine} 10: 73-78.

\bibitem[Rezaei et al. (2010)]{Rezaei} Rezaei, S., Tahmasebi, R., Mahmoodi, M., 2010. Estimation of $P(Y<X)$ for Generalized Pareto Distribution. {\it Journal of Statistical Planning and Inference} 140: 480-494.

\bibitem[Sara\c{c}o$\check{\mbox{g}}$lu and Kaya (2007)]{Saracoglu} Sara\c{c}o$\check{\mbox{g}}$lu, B., Kaya, M.F., 2007. Maximum likelihood estimation and confidence intervals of systemreliability for Gompertz distribution in stress-strength models. {\it Sel\c{c}uk Journal of Applied Mathematics} 8(2): 25-36.

\bibitem[Surles and Padgett (1998, 2001)]{SurlesPadgett} Surles, J.G., Padgett, W.J., 1998. Inference for $P(Y<X)$ in the Burr type $X$ model. {\it Journal of Applied Statistical Sciences} 7: 225-238.

\bibitem[Surles and Padgett (2001)]{SurlesPadgett2} Surles, J.G., Padgett, W.J., 2001. Inference for reliability and stress-strength for a scaled Burr-type $X$ distribution. {\it Life Time Data Analysis} 7: 187-200.

\bibitem[Tong (1977)]{Tong} Tong, H., 1977. On the estimation of $P(Y<X)$ for exponential families. {\it IEEE Transactions on Reliability} 26: 54-56.

\bibitem[Woodward and Kelley (1977)]{Woodward} Woodward, W.A., Kelley, G.D., 1977. Minimum variance unbiased estimation of $P(Y<X)$ in the normal case. {\it Technometrics} 19: 95-98.
\end{thebibliography}
\end{document}